\documentclass[11pt]{amsart}

\theoremstyle{plain}
\newtheorem{theorem}                {Theorem}      [section]
\newtheorem*{theorem*}                {Theorem \ref{thm:appl}}

\theoremstyle{definition}

\newtheorem{remark}       [theorem]  {Remark}
\newtheorem{definition}   [theorem]  {Definition}

\DeclareMathOperator{\trace}{trace}

 \DeclareMathOperator{\id}{I}

\DeclareMathOperator{\Span}{span}
 
\DeclareMathOperator{\cst}{constant}
 
 \DeclareMathOperator{\im}{Im}

\numberwithin{equation}{section}

\begin{document}

\title[Helix surfaces with parallel mean curvature]
{A classification result for helix surfaces with parallel mean curvature in product spaces}

\author{Dorel~Fetcu}

\address{Department of Mathematics and Informatics\\
Gh. Asachi Technical University of Iasi\\
Bd. Carol I no. 11 \\
700506 Iasi, Romania} \email{dfetcu@math.tuiasi.ro}

\curraddr{Department of Mathematics, Federal University of Bahia, Av.
Adhemar de Barros s/n, 40170--110 Salvador, BA, Brazil}

\thanks{The author was supported by a fellowship, BJT 373672/2013--6, offered by CNPq, Brazil.}



\begin{abstract} We determine all helix surfaces with parallel mean curvature vector field, which are not minimal or pseudo-umbilical, in spaces of type $M^n(c)\times\mathbb{R}$, where $M^n(c)$ is a simply-connected $n$-dimensional manifold with
constant sectional curvature $c$.
\end{abstract}

\maketitle

\section{Preliminaries and the main result}

Let us consider $M^n(c)$ a space form, i.e., a simply-connected $n$-dimensional manifold with
constant sectional curvature $c$, and the product manifold
$\bar M=M^n(c)\times\mathbb{R}$. 

We can deduce the expression of the curvature
tensor $\bar R$ of $\bar M$, from
$$
\langle\bar R(X,Y)Z,W\rangle=c\{\langle d\pi Y, d\pi Z\rangle\langle d\pi X, d\pi W\rangle-\langle d\pi X, d\pi Z\rangle\langle d\pi Y, d\pi W\rangle\},
$$
where $\pi:\bar M=M^n(c)\times\mathbb{R}\rightarrow M^n(c)$ is the projection map. We obtain
\begin{align}\label{eq:barR}
\bar R(X,Y)Z=&c\{\langle Y, Z\rangle X-\langle X, Z\rangle Y-\langle Y,\xi\rangle\langle Z,\xi\rangle X+\langle X,\xi\rangle\langle Z,\xi\rangle Y\\\nonumber &+\langle X,Z\rangle\langle Y,\xi\rangle\xi-\langle Y,Z\rangle\langle X,\xi\rangle\xi\},
\end{align}
where $\xi$ is the unit vector field tangent to $\mathbb{R}$.

Now, let $\Sigma^2$ be an isometrically immersed surface in $\bar M=M^n(c)\times\mathbb{R}$. The second fundamental form $\sigma$ of the surface is defined by the equation of Gauss
$$
\bar\nabla_XY=\nabla_XY+\sigma(X,Y),
$$
for any tangent vector fields $X$ and $Y$, where $\bar\nabla$ and $\nabla$ are the Levi-Civita connections on $\bar M$ and $\Sigma^2$, respectively. Then the mean curvature vector field $H$ of $\Sigma^2$ is given by $H=(1/2)\trace\sigma$. The shape operator $A$ and the normal connection $\nabla^{\perp}$ are defined by the equation of Weingarten
$$
\bar\nabla_XV=-A_VX+\nabla^{\perp}_XV,
$$
for any tangent vector field $X$ and any normal vector field $V$. 

We will also use the Gauss equation of the surface
\begin{align}\label{Gauss}
\langle R(X,Y)Z,W\rangle=&\langle\bar R(X,Y)Z,W\rangle+\langle\sigma(Y,Z),\sigma(X,W)\rangle\\\nonumber&-\langle\sigma(X,Z),\sigma(Y,W)\rangle,
\end{align}
and the Codazzi equation
\begin{equation}\label{Codazzi}
(\bar R(X,Y)Z)^{\perp}=(\nabla^{\perp}_X\sigma)(Y,Z)-(\nabla^{\perp}_Y\sigma)(X,Z),
\end{equation}
where $X$, $Y$, $Z$ and $W$ are tangent vector fields and $R$ is the curvature tensor corresponding to $\nabla$. 

\begin{definition} A surface $\Sigma^2$ in $M^n(c)\times\mathbb{R}$ is called a \textit{vertical cylinder} over $\gamma$ if $\Sigma^2=\pi^{-1}(\gamma)$, where $\pi:M^n(c)\times\mathbb{R}\rightarrow M^n(c)$ is the projection map and $\gamma:I\subset\mathbb{R}\rightarrow M^n(c)$ is a curve in $M^n(c)$.
\end{definition}

It is easy to see that vertical cylinders $\Sigma^2=\pi^{-1}(\gamma)$ are characterized by the fact that $\xi$ is tangent to $\Sigma^2$.

\begin{definition} A surface $\Sigma^2$ is called a \textit{helix surface} (or a \textit{constant angle surface}) if the angle function $\theta\in[0,\pi)$ between its tangent spaces and the vector field $\xi$ is constant on $\Sigma^2$. 
\end{definition}

A helix surface is characterized by the fact that the tangent part $T$ of $\xi$ has constant length.

\begin{definition} If the mean curvature vector $H$ of the surface $\Sigma^2$ is
parallel in the normal bundle, i.e.,
$\nabla^{\perp} H=0$, then $\Sigma^2$ is called a \textit{pmc surface}.
\end{definition}

\begin{definition} Let $\gamma:I\subset\mathbb{R}\rightarrow \bar M$ be a curve parametrized by
arc-length. Then $\gamma$ is called a {\it Frenet curve of
osculating order} $r$, $1\leq r\leq 2n$, if there exist $r$
orthonormal vector fields $\{X_1=\gamma',\ldots,X_r\}$ along
$\gamma$ such that
$$
\bar\nabla_{X_{1}}X_{1}=\kappa_{1}X_{2},\quad
\bar\nabla_{X_{1}}X_{i}=-\kappa_{i-1}X_{i-1} + \kappa_{i}X_{i+1},\quad\ldots\quad,
\bar\nabla_{X_{1}}X_{r}=-\kappa_{r-1}X_{r-1},
$$
for all $i\in\{2,\ldots,r-1\}$, where $\{\kappa_{1},\kappa_{2},\ldots,\kappa_{r-1}\}$ are positive
functions on $I$ called the {\it curvatures} of $\gamma$. A Frenet curve of osculating order $r$ is called a {\it helix of
order $r$} if $\kappa_i=\cst>0$ for $1\leq i\leq r-1$. A helix
of order $2$ is called a {\it circle}, and a helix of order $3$ is
simply called {\it helix}.
\end{definition}

When the ambient space is $M^2(c)\times\mathbb{R}$, a pmc surface $\Sigma^2$ is a surface with constant mean curvature (a \textit{cmc surface}). In order to study such surfaces, U.~Abresch and H.~Rosenberg introduced a holomorphic differential, now called the \textit{Abresch-Rosenberg differential}, and determined all complete cmc surfaces in $M^2(c)\times\mathbb{R}$ on which it vanishes (see \cite{AR}). They proved that there are four classes of such surfaces, denoted by $S_H^2$, $D_H^2$, $C_H^2$, and $P_H^2$, all of them described in detail in \cite{AR}. This holomorphic differential is the $(2,0)$-part of a quadratic form $Q$ defined on a cmc surface by
$$
Q(X,Y)=2\langle\sigma(X,Y),H\rangle-c\langle X,\xi\rangle\langle Y,\xi\rangle.
$$

Complete cmc helix surfaces in $M^2(c)\times\mathbb{R}$, and actually in all $3$-dimensional homogeneous spaces, were determined in
\cite[Theorem~2.2]{ER}, while \cite[Theorem~1]{CCCD} shows that there are no non-minimal pmc helix surfaces in $\mathbb{S}^3(1)\times\mathbb{R}$ with $0<|T|<1$.
 
In our paper, we consider non-minimal pmc helix surfaces in $M^n(c)\times\mathbb{R}$ and prove the following classification theorem.

\begin{theorem}\label{main} Let $\Sigma^2$ be a non-minimal pmc helix surface with Gaussian curvature $K$ and mean curvature vector field $H$ in $\bar M=M^n(c)\times\mathbb{R}$, with $c\neq 0$, and let $T$ be the tangent part of the unit vector field $\xi$ tangent to $\mathbb{R}$. Then one of the following holds$:$
\begin{enumerate}

\item $\Sigma^2$ is a minimal surface in a non-minimal totally umbilical hypersurface of $M^n(c)$$;$

\item $\Sigma^2$ is a surface with constant mean curvature in a $3$-dimensional totally umbilical or totally geodesic submanifold of $M^n(c)$$;$

\item $\Sigma^2$ is a vertical cylinder over a circle in $M^n(c)$ with curvature $\kappa=2|H|$$;$

\item $c<0$ and $\Sigma^2$ lies in $M^2(c)\times\mathbb{R}$. Moreover, $4|H|^2+c|T|^2=0$, $K=c(1-|T|^2)<0$, and the Abresch-Rosenberg differential vanishes on $\Sigma^2$. If the surface is complete, then it is a cmc surface of type $P_H^2$ in $M^2(c)\times\mathbb{R}$$;$

\item $c>0$ and locally $\Sigma^2$ is the standard product $\gamma_1\times\gamma_2$, where $\gamma_1:I\subset\mathbb{R}\rightarrow M^4(c)\times\mathbb{R}$ is a helix in $M^4(c)\times\mathbb{R}$ with curvatures $\kappa_1=\sqrt{c(1-|T|^2)}$ and $\kappa_2=|T|\sqrt{c}$, and $\gamma_2:I\subset\mathbb{R}\rightarrow M^4(c)\subset M^4(c)\times\mathbb{R}$ is a circle in $M^4(c)$ with curvature $\kappa=\sqrt{4|H|^2+c(1-|T|^2)}$. If the surface is complete, then the above decomposition holds globally$.$

\end{enumerate}
\end{theorem}

\noindent {\bf Acknowledgments.} The author would like to thank the Department of Mathematics of the Federal University of Bahia in Salvador for providing a very stimulative work environment during the preparation of this paper.

\section{The proof of Theorem \ref{main}}

The map $p\in\Sigma^2\rightarrow(A_H-\mu\id)(p)$ is analytic, and,
therefore, either $\Sigma^2$ is a pseudo-umbilical surface, i.e., $A_H=|H|^2\id$,
or $H$ is an umbilical direction on a closed set without interior
points. We shall denote by $W$ the set of points where $H$ is not an
umbilical direction, which, in the second case, is an open dense set in
$\Sigma^2$.

If $\Sigma^2$ is a pmc surface in $\mathbb{M}^n(c)\times\mathbb{R}$, with mean curvature vector field $H$, then either $\Sigma^2$ is
pseudo-umbilical, i.e., $H$ is an umbilical direction everywhere, or, at any point in $W$, there exists a local orthonormal frame field that diagonalizes $A_U$ for
any normal vector field $U$ defined on $W$ (see \cite[Lemma 1]{AdCT}). If $\Sigma^2$ is a pseudo-umbilical pmc surface in $\mathbb{M}^n(c)\times\mathbb{R}$, then it lies in $\mathbb{M}^n(c)$, i.e. $|T|=0$ (see \cite[Lemma 3]{AdCT}).

From \cite[Remark~1]{AdCT}, we also know that, since $H$ is parallel, the immersion of $\Sigma^2$ in $M^n(c)\times\mathbb{R}$ is analytic, i.e., the functions of two variables that locally define the immersion are real analytic. Therefore, it satisfies a principle of unique continuation and, as a consequence, it cannot vanish on an open connected subset of $\Sigma^2$ unless it vanishes identically.

Now, when $|T|=0$ our surface lies in $M^n(c)$ and we obtain the first two items of the theorem using \cite[Theorem~4]{Y}. 

When $|T|=1$ the vector field $\xi$ is tangent to the surface and this means that $\Sigma^2$ is a pmc vertical cylinder over a curve $\gamma$ in $M^n(c)$. Since $\xi$ is parallel in $M^n(c)\times\mathbb{R}$, it follows that $\sigma(\xi,\xi)=0$ and then we easily get that $\gamma$ is a circle in $M^n(c)$ with curvature $\kappa=2|H|$.

Henceforth, let us assume that $0<|T|<1$. It follows that $H$ is not an umbilical direction on an open dense set $W$. We will work on this set and then we will extend our results throughout $\Sigma^2$ by continuity. 

Consider a global orthonormal frame field $\{E_1=T/|T|,E_2\}$ on $\Sigma^2$, and let $N$ be the normal part of $\xi$. Then, since $\Sigma^2$ is a helix surface, it follows that $\nabla_{E_1}E_1=\nabla_{E_1}E_2=0$ and, as $\bar\nabla_X\xi=0$ implies $\nabla_XT=A_NX$ and $\sigma(T,X)=-\nabla^{\perp}_XN$, that $A_NE_1=0$ (see \cite[Proposition~2.1]{RH}). We also have
\begin{align*}
\langle\nabla_{E_2}E_2,E_1\rangle &=-\langle E_2,\nabla_{E_2}E_1\rangle=-\frac{1}{|T|}\langle E_2,A_NE_2\rangle=-\frac{1}{|T|}\langle\sigma(E_2,E_2),N\rangle\\ &=-\frac{1}{|T|}\langle 2H-\sigma(E_1,E_1),N\rangle=-\frac{2\langle H,N\rangle}{|T|},
\end{align*}
which means that 
\begin{equation}\label{nabla}
\nabla_{E_2}E_2=-\frac{2\langle H,N\rangle}{|T|}E_1\quad\textnormal{and}\quad\nabla_{E_2}E_1=\frac{2\langle H,N\rangle}{|T|}E_2.
\end{equation}

Since $A_NE_1=0$, we also get that $A_NE_2=2\langle H,N\rangle$. From the Ricci equation
$$
\langle R^{\perp}(X,Y)U,V\rangle=\langle [A_U,A_V]X,Y\rangle+\langle\bar R(X,Y)U,V\rangle,
$$
where $X$ and $Y$ are tangent vector fields and $U$ and $V$ are normal vector fields, we obtain $[A_H,A_N]=0$, which means that $\langle H,N\rangle\langle A_HE_1,E_2\rangle=0$. Since we also have 
$$
E_2(\langle H,N\rangle)=\langle H,\nabla^{\perp}_{E_2}N\rangle=-|T|\langle H,\sigma(E_1,E_2)\rangle=-|T|\langle A_HE_1,E_2\rangle,
$$
one sees that $\langle A_HE_1,E_2\rangle=0$. Moreover, again using the Ricci equation, we have $[A_H,A_U]=0$ for any normal vector field $U$ and then, since $H$ is not umbilical, $\sigma(E_1,E_2)=0$ and $\nabla^{\perp}_{E_2}N=0$. 

Next, let us denote $\lambda_i=\langle A_HE_i,E_i\rangle$, $i\in\{1,2\}$, the eigenfunctions of $A_H$ and in the following we will compute $E_j(\lambda_i)$, $i,j\in\{1,2\}$.

Using \eqref{eq:barR}, \eqref{nabla} and $\sigma(E_1,E_2)=0$, we get
$$
\nabla^{\perp}_{E_2}\sigma(E_1,E_1)=(\bar R(E_2,E_1)E_1)^{\perp}+2\sigma(E_1,\nabla_{E_2}E_1)=0
$$
and then, from the Codazzi equation \eqref{Codazzi}, 
\begin{align}\label{eq:e1sigma}
\nabla^{\perp}_{E_1}\sigma(E_2,E_2)=&(\bar R(E_1,E_2)E_2)^{\perp}+2\sigma(E_2,\nabla_{E_1}E_2)+\nabla^{\perp}_{E_2}\sigma(E_1,E_1)\\\nonumber&-\sigma(\nabla_{E_2}E_1,E_2)-\sigma(E_1,\nabla_{E_2}E_2)\\\nonumber=&-c|T|N+\frac{2\langle H,N\rangle}{|T|}(\sigma(E_1,E_1)-\sigma(E_2,E_2)).
\end{align}
Therefore, since $H$ is parallel, we have
\begin{equation}\label{eq:e1lambda}
E_1(\lambda_1)=-E_1(\lambda_2)=-\langle\nabla^{\perp}_{E_1}\sigma(E_2,E_2),H\rangle=\frac{\langle H,N\rangle}{|T|}(4|H|^2+c|T|^2-4\lambda_1)
\end{equation}
and
\begin{equation}\label{eq:e2lambda}
E_2(\lambda_1)=-E_2(\lambda_2)=\langle\nabla^{\perp}_{E_2}\sigma(E_1,E_1),H\rangle=0.
\end{equation}

From \cite[Proposition~1.4]{FR}, we know that
$$
\frac{1}{2}\Delta|T|^2=|A_N|^2+K|T|^2-2\langle A_HT,T\rangle
$$
and then, since $|T|=\cst\neq 0$, the Gaussian curvature $K$ of $\Sigma^2$ is given by
\begin{equation}\label{eq:K}
K=2\lambda_1-\frac{4\langle H,N\rangle^2}{|T|^2}.
\end{equation}
Since $\nabla^{\perp}_XN=-\sigma(X,T)$ implies that
\begin{equation}\label{eq:eHN}
E_1(\langle H,N\rangle)=-|T|\lambda_1\quad\textnormal{and}\quad E_2(\langle H,N\rangle)=0,
\end{equation}
from \eqref{eq:e1lambda} and \eqref{eq:e2lambda}, we obtain
\begin{equation}\label{eq:eK}
E_1(K)=\frac{2\langle H,N\rangle}{|T|^2}(4|H|^2+c|T|^2)\quad\textnormal{and}\quad E_2(K)=0.
\end{equation}

The fact that $\Sigma^2$ is not pseudo-umbilical implies that it lies in $M^4(c)\times\mathbb{R}$ (see \cite[Theorem~1]{AdCT}). 

In order to describe our surface by taking advantage of the above formulas, we will first consider the case when $H\parallel N$ on an open connected subset $W_0$ of $\Sigma^2$. This means that $\langle H,N\rangle=\pm|H||N|=\cst$ and, from \eqref{eq:eHN} one obtains that $\lambda_1=0$. Then, from \eqref{eq:e1lambda} and \eqref{eq:K}, we get $4|H|^2+c|T|^2=0$ and $K=c(1-|T|^2)$ on $W_0$. 

Let $\{E_3=H/|H|,E_4,E_5\}$ be a global orthonormal frame field in the normal bundle. Then, since $\sigma(E_1,E_2)=0$, on $W_0$ we have 
$$
A_3=\left(\begin{array}{cc}0&0\\0&2|H|\end{array}\right),\quad A_4=\left(\begin{array}{cc}a&0\\0&-a\end{array}\right),\quad A_5=\left(\begin{array}{cc}b&0\\0&-b\end{array}\right),
$$
where $A_{\alpha}=A_{E_{\alpha}}$, and, from the Gauss equation \eqref{Gauss}, 
$$
K=c(1-|T|^2)-a^2-b^2.
$$
But, as we have seen, $K=c(1-|T|^2)$ on $W_0$, which implies that $a=b=0$ on $W_0$. Consider the subbundle $L=\Span\{\im\sigma\}=\Span\{H\}$ in the normal bundle. It follows that $\xi\in T\Sigma\oplus L$ and $T\Sigma^2\oplus L$ is parallel with respect to $\bar\nabla$ and invariant by $\bar R$. Therefore, we use \cite[Theorem~2]{ET} to show that $W_0$ lies in $M^2(c)\times\mathbb{R}$. From the analyticity of the immersion of $\Sigma^2$ in $M^4(c)\times\mathbb{R}$, it follows that $\Sigma^2$ lies in $M^2(c)\times\mathbb{R}$ (see \cite[Remark~1]{AdCT} for more details). As a direct consequence, we have $H\parallel N$ on $\Sigma^2$ and then $\lambda_1=0$, $4|H|^2+c|T|^2=0$ and $K=c(1-|T|^2)<0$ on $\Sigma^2$. Now, it is easy to verify that the Abresch-Rosenberg differential vanishes on the surface. Moreover, if $\Sigma^2$ is complete, all these properties of $\Sigma^2$ lead to the conclusion that our surface is a cmc surface of type $P_H^2$. 

Next, we will consider the remaining case, when $H\parallel N$ only at isolated points. We can then define a local orthonormal frame field 
$$
\Big\{E_3=\frac{1}{|H|\sin\beta}H-\frac{\cot\beta}{|N|}N,E_4=\frac{N}{|N|},E_5\Big\}
$$
in the normal bundle, where $\beta\in(0,\pi)$ is the angle between $H$ and $N$. With respect to $\{E_1,E_2\}$ we can write
$$
A_3=\left(\begin{array}{cc}\frac{\lambda_1}{|H|\sin\beta}&0\\0&2|H|\sin\beta-\frac{\lambda_1}{|H|\sin\beta}\end{array}\right), A_4=\left(\begin{array}{cc}0&0\\0&\frac{2\langle H,N\rangle}{|N|}\end{array}\right),A_5=\left(\begin{array}{cc}\lambda&0\\0&-\lambda\end{array}\right)
$$
and then, from the Gauss equation \eqref{Gauss}, one obtains
\begin{equation}\label{eq:Klambda}
K=c(1-|T|^2)+2\lambda_1-\frac{\lambda_1^2}{|H|^2\sin^2\beta}-\lambda^2.
\end{equation}
Using equation \eqref{eq:K}, we have
\begin{equation}\label{eq:1}
\lambda^2=c(1-|T|^2)-\frac{\lambda_1^2}{|H|^2\sin^2\beta}+\frac{4\langle H,N\rangle^2}{|T|^2}
\end{equation}
and 
\begin{equation}\label{eq:e11}
2\lambda E_1(\lambda)=E_1\Big(-\frac{\lambda_1^2}{|H|^2\sin^2\beta}+\frac{4\langle H,N\rangle^2}{|T|^2}\Big).
\end{equation}

Next, we shall compute $E_1(\lambda)$. Using equation \eqref{eq:e1sigma}, the fact that $H$ is parallel and $\nabla^{\perp}_XN=-\sigma(X,T)$, we obtain
\begin{align*}
E_1(\lambda)&=-E_1(\langle\sigma(E_2,E_2),E_5\rangle)=-\langle\nabla^{\perp}_{E_1}\sigma(E_2,E_2),E_5\rangle-\langle\sigma(E_2,E_2),\nabla^{\perp}_{E_1}E_5\rangle\\&=-\frac{4\langle H,N\rangle}{|T|}\lambda-\langle\sigma(E_2,E_2),E_3\rangle\langle\nabla^{\perp}_{E_1}E_5,E_3\rangle-\langle\sigma(E_2,E_2),E_4\rangle\langle\nabla^{\perp}_{E_1}E_5,E_4\rangle\\&=-\frac{4\langle H,N\rangle}{|T|}\lambda+\Big(2|H|\sin\beta-\frac{\lambda_1}{|H|\sin\beta}\Big)\langle E_5,\nabla^{\perp}_{E_1}E_3\rangle+\frac{2\langle H,N\rangle}{|N|}\langle E_5,\nabla^{\perp}_{E_1}E_4\rangle\\&=-\frac{4\langle H,N\rangle}{|T|}\lambda-\frac{|T|}{|H||N|}\frac{\cos\beta}{\sin^2\beta}\lambda\lambda_1.
\end{align*} 

Replacing in \eqref{eq:e11} and using \eqref{eq:e1lambda} and \eqref{eq:eHN}, we get, after a straightforward computation
\begin{equation}\label{eq:e11final}
\langle H,N\rangle\Big(2\lambda_1-\frac{\lambda_1^2}{|H|^2\sin^2\beta}-\lambda^2\Big)=0,
\end{equation}
which, together with \eqref{eq:Klambda}, leads to $\langle H,N\rangle(K-c(1-|T|^2))=0$,
that, taking \eqref{eq:eK} into account, can be written as $E_1((K-c(1-|T|^2))^2)=0$.
Again from \eqref{eq:eK}, we have $E_2((K-c(1-|T|^2))^2)=0$. It follows that $K-c(1-|T|^2)=\cst$ 
and then, using \eqref{eq:K}, one obtains
$$
E_1\Big(2\lambda_1-\frac{4\langle H,N\rangle^2}{|T|^2}\Big)=0.
$$

Now, from equations \eqref{eq:e1lambda} and \eqref{eq:eHN}, we have $(4|H|^2+c|T|^2)\langle H,N\rangle=0$. We will consider two cases as $4|H|^2+c|T|^2=0$ or $4|H|^2+c|T|^2\neq 0$.

\textbf{Case I: $4|H|^2+c|T|^2=0$.} Let us assume that $\langle H,N\rangle=0$ on an open connected set $W_0$. From \eqref{eq:eHN} it follows that $\lambda_1=0$ and then, from \eqref{eq:K}, that $K=0$ on $W_0$. But, from \eqref{eq:Klambda}, we have $K=c(1-|T|^2)-\lambda^2$, which means that $\lambda^2=c(1-|T|^2)$. This implies $c>0$ and that is a contradiction, since $4|H|^2+c|T|^2=0$. Therefore $\langle H,N\rangle=0$ on a closed set without interior points and then, from \eqref{eq:e11final}, we have 
$$
2\lambda_1-\frac{\lambda_1^2}{|H|^2\sin^2\beta}-\lambda^2=0
$$
and $K=c(1-|T|^2)$ on an open dense set. Since $4|H|^2+c|T|^2=0$, from \eqref{eq:K}, one obtains $\sin^2\beta=2\lambda_1/(c|N|^2)$ and then $\lambda^2=2\lambda_1/|T|^2$, which means that $\lambda_1=0$, $\lambda=0$ and $\sin^2\beta=0$ on an open dense set. the last identity shows that $H\parallel N$ on an open dense set, which is a contradiction.

\textbf{Case II: $4|H|^2+c|T|^2\neq 0$.} In this case, $\langle H,N\rangle=0$ on an open dense set and, from \eqref{eq:eHN}, it follows that $\lambda_1=0$ and then, from \eqref{eq:K}, we have $K=0$ and, from \eqref{eq:1}, $\lambda^2=c(1-|T|^2)$, which implies $c>0$. The shape operator is given by
$$
A_3=\left(\begin{array}{cc}0&0\\0&2|H|\end{array}\right),\quad A_4=\left(\begin{array}{cc}0&0\\0&0\end{array}\right),\quad A_5=\left(\begin{array}{cc}\lambda&0\\0&-\lambda\end{array}\right).
$$
From the equations \eqref{nabla}, we see that $\nabla E_1=\nabla E_2=0$ and we then can apply the de Rham Decomposition Theorem (\cite{KN}) to show that locally $\Sigma^2$ is the standard product $\gamma_1\times\gamma_2$ of two Frenet curves $\gamma_1:I\subset\mathbb{R}\rightarrow M^4(c)\times\mathbb{R}$ and $\gamma_2:I\subset\mathbb{R}\rightarrow M^4(c)\times\mathbb{R}$ such that $\gamma'_1=E_1$ and $\gamma_2'=E_2$. We note that $\gamma_2$ actually lies in $M^4(c)$. If the surface is complete, then this decomposition holds globally.

Next, we will characterize $\gamma_1$ and $\gamma_2$ by using their Frenet equations. Let $\{X_1^1=E_1,X_2^1,\ldots,X_r^1\}$ be the Frenet frame field of $\gamma_1$. We have that
$$
\bar\nabla_{E_1}E_1=\sigma(E_1,E_1)=\lambda E_5,
$$
and then, from the first Frenet equation, it follows that $\kappa_1=|\lambda|=\sqrt{c(1-|T|^2)}$ and $X_2^1=\pm E_5$.

Next, we have $\langle\nabla^{\perp}_{E_1}E_5,E_3\rangle=0$, since $E_3=H/|H|$ is parallel, and $\langle\nabla^{\perp}_{E_1}E_5,E_4\rangle=-\langle E_5,\nabla^{\perp}_{E_1}E_4\rangle=(|T|/|N|)\langle E_5,\sigma(E_1,E_1)\rangle=(|T|/|N|)\lambda$. Therefore, one obtains
$$
\bar\nabla_{E_1}X_2^1=\pm\bar\nabla_{E_1}E_5=\mp A_5E_1\pm\nabla^{\perp}_{E_1}E_5=\mp\lambda E_1\pm\frac{|T|}{|N|}\lambda E_4.
$$
From the second Frenet equation, we get $\kappa_2=(|T|/|N|)|\lambda|=|T|\sqrt{c}$ and $X_3^1=\pm E_4$.

It follows that
$$
\bar\nabla_{E_1}X_3^1=\pm\bar\nabla_{E_1}E_4=\pm\nabla^{\perp}_{E_1}E_4=\mp(|T|/|N|)\sigma(E_1,E_1)=\mp(|T|/|N|)\lambda E_5=-\kappa_2X_2^1
$$
and we have just proved that $\gamma_1$ is a helix.

Now, let $\{X_1^2=E_2,X_2^2,\ldots,X_r^2\}$ be the Frenet frame field of $\gamma_2$. Then, from
$$
\bar\nabla_{E_2}E_2=\sigma(E_2,E_2)=2|H|E_3-\lambda E_5,
$$
and the first Frenet equation of $\gamma_2$, one obtains 
$$
\kappa=\sqrt{4|H|^2+\lambda^2}=\sqrt{4|H|^2+c(1-|T|^2)}
$$ 
and $X_2^2=(1/\kappa)(2|H|E_3-\lambda E_5)$.

It is easy to verify, using $\bar\nabla\xi=0$, that $\nabla^{\perp}_{E_2}E_5=0$. Then, since $E_3$ is parallel, we have $\bar\nabla_{E_2}X_2^2=-(1/\kappa)(2|H|A_3E_2+\lambda A_3E_2)=-\kappa E_2$ and we conclude.

\begin{remark} From the proof of Theorem \ref{main} it is easy to see that the only non-minimal pmc helix surfaces in $M^3(c)\times\mathbb{R}$ are those given by the first four cases of our theorem.
\end{remark}

\end{document}